# A ROBUST OPTIMIZATION APPROACH MODEL FOR A MULTI-VACCINE MULTI-ECHELON SUPPLY CHAIN

A. BOUCHENINE[1,2*], I. ALMARAJ[1,2]

[1]Department of Industrial and Systems Engineering King Fahd University of Petroleum and Minerals, 5063, Dhahran 31261, Saudi Arabia
[2]Interdisciplinary Research Center for Smart Mobility and Logistics, King Fahd University of Petroleum and Minerals, 5063, Dhahran 31261, Saudi Arabia
b92.abderrahmen@gmail.com, g202113210@kfupm.edu.sa,
almaraj@kfupm.edu.sa

### ABSTRACT

This research investigates a multi-product, multi-echelon and multi-period vaccine supply chain (SC) network model under uncertainty and quality inspection errors. The objective function seeks optimizing the total cost (economic cost) of the SC. Moreover, the proposed model is formulated as a Mixed Integer Linear Programming problem (MILP) under multiple sources of uncertain parameters including demand, inspection errors, vaccine waste generated in healthcare centers, and defective treatment rate of vaccine waste. To provide meaningful solutions which are robust against the future fluctuation of parameters, the robust optimization approach is utilized to incorporate the decision maker's risk attitude under different type of uncertainty sets namely, "box", "polyhedral" and combination of "interval polyhedral". The performance of the proposed model is demonstrated through an illustrative example. The results show the effect of different types of uncertainties on the overall objective function. Managerial insights and research implications in terms of reverse vaccine SC is advised and future research directions are proposed.

## 1  INTRODUCTION AND LITERATURE REVIEW

Vaccines are one of the most cost-effective methods to contain outbreaks. Over the years, vaccination has proven to be the most efficient way to prevent and control the spread of infectious diseases, Yang et al.[1]. Vaccine supply chain has drastically affected by the emergence of COVID-19 pandemic in the last two-years. Therefore, distribution of vaccines become a challenge for policymakers. As a result, the World Health Organization (WHO) defined priorities for containing epidemics by distributing mass vaccination (especially for developing countries) into three main priorities namely, products and packaging, immunization supply system efficiency, and environmental impact of immunization supply systems, Duijzer et al.[2]. Recently, many researchers have extensively addressed these preferences in the literature in order to tackle the vaccine supply chain problems through its life cycle. Most importantly, a focus has been shaded on four components of the vaccine supply chain; product, production, allocation and distribution, Duijzer et al.[2]. Vaccine delivery may take different levels to the patients, for instance, from supplier to distribution centers and finally to medical centers (hospitals, school/university medical centers).Thus, mass vaccination planning problem is a challenge for both researchers and practitioners, Tang et al.[3]. Another issue that may arise during a pandemic containment is the huge amounts of medical waste generated from healthcare centers due to patient's vaccination campaign. This kind of waste is hazardous and should be treated accordingly, Dasaklis et al.[4]. Therefore, vaccine waste is another threat that must be treated carefully and get disposed in such a way that they do not pose a risk for the medical personnel and people engaged in the containment effort, Dasaklis et al.[4]. Vaccine supply chain is prone to disruption by many factors, one of the major causes that may lead to shortages and stockouts of vaccines in the market is imperfect or contaminated items (defective items). As results, vaccine inspection could take place in different stages along its distribution to the end user.

---

* Corresponding Author: b92.abderrahmen@gmail.com, g202113210@kfupm.edu.sa (A.Bouchenine).

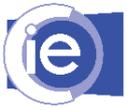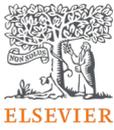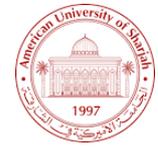

Since the appearance of COVID-19 outbreak, many research have been conducted on vaccine supply chain. For instance, Gilani and Sahebi [5], presented a mathematical model which propose a robust data-driven model based on a polyhedral uncertainty set to tackle the vaccine distribution as an uncertain parameter. Manupati et al.[6], suggested a decision support in order to distribute the vaccine and give access to it. In terms of vaccine distribution, many aspects have been investigated in the literature; Bertsimas et al.[7],tackled the question of 'Where to locate COVID-19 mass vaccination facilities' with a novel data-driven approach to optimize COVID-19 vaccine distribution. Duijzer et al.[2], proposed a classification for the literature on vaccine logistics to structure new field and identify promising research directions.

The existing literature on vaccine quality control along with its waste management need to be extended, since quite few works have combined these issues together in the literature. In contrast, cold supply chain has been drastically carried out during the past years in different setting.

To this end, the vaccine SC with inspection errors along with its waste are critical aspect for stakeholders during an outbreak like COVID-19. Which may limit the access to the vaccine items and cause shortages at the customer level (healthcare centers). According to the conducted literature, vaccine SC with inspection errors, and its waste management have not been addressed before in the literature, which motivated us to investigate the proposed problem.

The remainder of the manuscript is organized as follows; we first define the problem and formulate the deterministic model, then the robust convex optimization approach is introduced and the deterministic model is reformulated as a robust counterpart program for the three uncertainty sets, box, polyhedral and combined interval-polyhedral. Next, to validate the proposed model, an illustrative example is introduced; computational results and analysis are conducted. Finally, the paper is concluded with comprehensive discussion along with conclusion and future extensions.

## 2    PROBLEM STATEMENT AND FORMULATION

In this study, we consider a vaccine supply chain system consisting of multiple periods, multiproduct, and echelons. The flow of materials depicted in figure-1, and can be described as follows; the vaccine orders are received from suppliers(s) to a set of healthcare distribution centers for inspection, and then are moved to the healthcare centers according to the demand. Furthermore, due to a massive vaccination campaign, a considerable amount of vaccine waste is generated. Consequently, the vaccine wastes resulted from healthcare centers are shipped to vaccine waste storage centers along with defective items from the distribution centers. It is assumed that vaccine storage centers have an inventory in order to handle the massive vaccine waste generated from healthcare centers. Another reason of contemplating the inventory in the storage centers is that vaccine waste treatment centers have a defective treatment rate, which cannot treat all the generated vaccine waste directly due to machine breakdowns, periodic maintenance.

Next, vaccine wastes are transported to potential treatment centers, these facilities in practice are dedicated to treat special medical waste generated from hospitals and clinics with special technology. Finally, the treated wastes are disposed through landfills. The supply chain facilities are considered to be fixed and predefined. The provided healthcare distribution centers serve as inspection unit in which the vaccine orders received are inspected and screened for any defects, in addition, type I & II errors ($\gamma_1$ & $\gamma_2$) assumed to occur during the inspection and considered as uncertain parameters. The objective is to minimize the total cost of the supply chain network to better serve the healthcare centers with vaccine items and get dispose of vaccine waste generated in an efficient manner.

The main assumptions of the problem are as follows:
- It is assumed that items in the potential healthcare distribution centers are subject to partial disruption due to inspection errors.
- Healthcare distribution centers consists of inspection units and warehouses. The received vaccine items are all inspected and screened for any defects in the inspection units and then moved to the warehouses. While defective items are discarded directly through the storage waste centers.

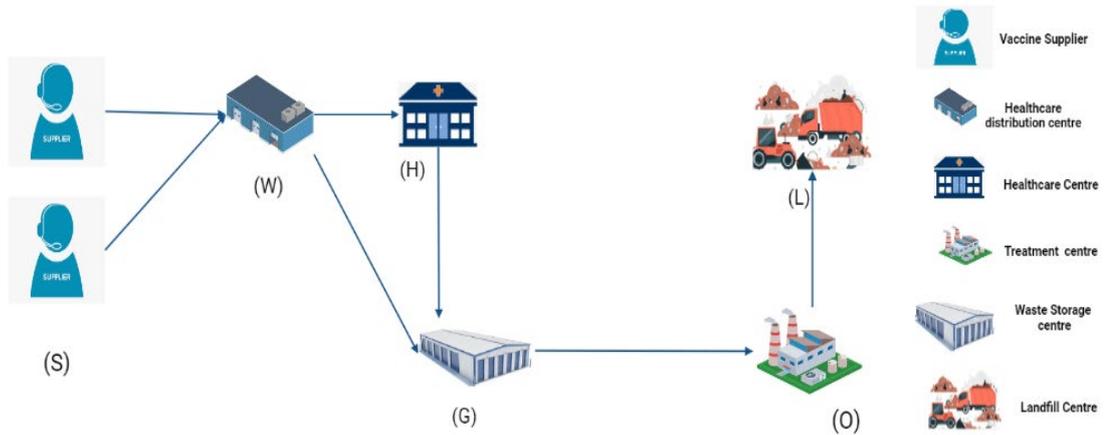

Figure 1: Supply chain network for the proposed problem.

## 2.1 Notation

| Sets | Parameters |
|---|---|
| S  Set of possible suppliers $s \in S$. | $FS_{st}$  The fixed cost of selecting supplier $s \in S$ at time t. |
| W Set of potential distribution centers $w \in W$. | $FW_{wt}$  The fixed cost of opening a warehouse $w \in W$ at time t. |
| H Set of healthcare centers $h \in H$. | |
| O Set of potential treatment centers $o \in O$. | $FG_{gt}$  The fixed cost of opening vaccine waste storage center $g \in G$ at time t. |
| G Set of potential storage centers $g \in G$. | |
| L Set of potential landfills $l \in L$. | $FO_{ot}$  The fixed cost of opening a treatment center $o \in O$ at time t. |
| P Set of vaccine products types $p \in P$. | |
| T Set of time periods $t \in T$. | |

| Decision variables | |
|---|---|
| $QS_{swpt}$ | Amount of vaccine ordered from supplier $s \in S$ to distribution centers $w \in W$ at time t. |
| $QW_{whpt}$ | Amount of vaccine shipped from distribution centers $w \in W$ to healthcare centers $h \in H$ at time t. |
| $QD_{wgpt}$ | Amount of defective vaccine shipped from distribution centers $w \in W$ to vaccine waste storage center $g \in G$ at time t. |
| $QU_{hpt}$ | Amount of unmet demand of vaccine items to the healthcare center $h \in H$ at time t. |
| $QH_{hgt}$ | Amount of vaccine waste shipped from healthcare centers $h \in H$ to vaccine waste storage centers $g \in G$ at time t. |
| $QG_{got}$ | Amount of vaccine waste shipped from waste storage centers $g \in G$ to treatment facility $o \in O$ at time t. |
| $QO_{olt}$ | Amount of treated vaccine waste shipped from treatment centers $o \in O$ to landfill $l \in L$ at time t. |
| $INV_{gt}$ | Inventory level of vaccine waste at the waste storage center $g \in G$ at time t. |
| $Y_{wt}$ | = 1: If an order is placed at time t. 0, otherwise. |
| $S_{st}$ | = 1: If a supplier $s \in S$ is selected at time t, 0 otherwise. |
| $W_{wt}$ | = 1: If a distribution center $w \in W$ is opened at time t, 0 otherwise. |
| $G_{gt}$ | = 1: If storage center is opened at time t, 0 otherwise. |
| $O_{ot}$ | = 1: If a treatment center $o \in O$ is opened at time t, 0 otherwise. |
| $L_{lt}$ | = 1: If landfill $l \in L$ is opened at time t, 0 otherwise. |

| Uncertain parameters | |
|---|---|
| $\widetilde{D}_{hpt}$ | The demand of vaccine items $p \in P$ in the healthcare center $h \in H$ at time t. |
| $\widetilde{VW}_{ht}$ | The amount of vaccine waste generated in the healthcare center $h \in H$ at time t. |
| $\widetilde{\Theta}_o$ | The defective treatment rate in the potential treatment center $o \in O$. |
| $\tilde{\gamma}_1$ | Probability of type I error of inspection. |
| $\tilde{\gamma}_2$ | Probability of type II error of inspection. |
| $FL_{lt}$ | The fixed cost of opening a landfill $l \in L$ at time t. |
| $D_{hpt}$ | The demand of vaccine items $p \in P$ in the healthcare center $h \in H$ at time t. |
| $VW_{ht}$ | The amount of vaccine waste generated in the healthcare center $h \in H$ at time t. |
| $KW_{wpt}$ | Ordering cost per lot size of vaccine items $p \in P$ incurred at the distribution center $w \in W$ at time t. |
| $HCW_{wp}$ | The inventory holding cost of vaccine items $p \in P$ at the distribution center $w \in W$. |
| $p_h$ | The shortage cost of vaccine items in the healthcare center $h \in H$. |
| $c_{sp}$ | The purchasing cost per vaccine items $p \in P$ at the supplier $s \in S$. |
| $HCG_g$ | The holding cost of vaccine waste in the waste storage centers $g \in G$. |
| $TC_o$ | Treatment cost incurred in treatment center $o \in O$. |
| $\Theta_o$ | The defective treatment rate in the potential treatment center $o \in O$. |
| $I_{wp}$ | Inspection cost of the vaccine items $p \in P$ incurred at the distribution center $w \in W$. |
| $p$ | The probability of defective items. |
| $(1 - p)$ | The probability of non-defective items. |
| $\bar{p}$ | The apparent probability of defective items. |
| $(1 - \bar{p})$ | The apparent probability of non-defective items. |
| $\gamma_1$ | Probability of type I error of inspection. |
| $\gamma_2$ | Probability of type II error of inspection. |
| $TS_{sw}$ | Transportation cost of vaccine items from supplier $s \in S$ to distribution center $w \in W$. |
| $TW_{wh}$ | Transportation cost of vaccine from distribution center $w \in W$ to healthcare center $h \in H$. |
| $TG_{wg}$ | Transportation cost of defective vaccine items from distribution center $w \in W$ to waste storage center $o \in O$. |
| $TH_{hg}$ | Transportation cost of vaccine waste from healthcare center $h \in H$ to waste storage center $g \in G$. |
| $TM_{go}$ | Transportation cost of vaccine waste from the waste storage center $g \in G$ to treatment center $o \in O$. |
| $TO_{ol}$ | Transportation cost of treated vaccine waste from treatment center $o \in O$ to landfill $l \in L$. |
| $CS_{sp}$ | The supplier capacity of vaccine $p \in P$. |
| $CW_{wp}$ | The distribution center capacity of vaccine $p \in P$. |
| $CG_g$ | The waste storage center capacity $g \in G$. |
| $CO_o$ | The treatment center capacity, $o \in O$. |
| $CL_l$ | The landfills capacity, $l \in L$. |
| $M$ | Big number. |

The deterministic problem is formulated as follows:

$$Minimize\ Z_1 = \sum_{t \in T}\sum_{s \in S} FS_{st} \times S_{st} + \sum_{t \in T}\sum_{w \in W} FW_{wt} \times W_{wt} + \sum_{t \in T}\sum_{g \in G} FG_{gt} \times G_{gt}$$
$$+ \sum_{t \in T}\sum_{o \in O} FO_{ot} \times O_{ot} + \sum_{t \in T}\sum_{l \in L} FL_{lt} \times L_{lt} + \sum_{t \in T}\sum_{s \in S}\sum_{w \in W}\sum_{p \in P} QS_{swpt} \times TS_{swp}$$
$$+ \sum_{t \in T}\sum_{w \in W}\sum_{p \in P} KW_{wpt} \times Y_{wt} + \sum_{t \in T}\sum_{s \in S}\sum_{w \in W}\sum_{p \in P} QS_{swpt} \times c_{spt}$$
$$+ \sum_{t \in T}\sum_{s \in S}\sum_{w \in W}\sum_{p \in P} (QS_{swpt} \times I_{wp}) + \sum_{t \in T}\sum_{s \in S}\sum_{w \in W}\sum_{p \in P} (QS_{swpt} \times HCW_w)$$

$$+ \sum_{t \in T} \sum_{w \in W} \sum_{h \in H} \sum_{p \in P} QW_{whpt} \times (1 - \bar{p}) \times TW_{whp} \quad (1)$$

$$+ \sum_{t \in T} \sum_{w \in W} \sum_{h \in H} \sum_{p \in P} QU_{hpt} \times p_h + \sum_{t \in T} \sum_{w \in W} \sum_{g \in G} \sum_{p \in P} QD_{wgpt} \times \bar{p} \times TG_{wg}$$

$$+ \sum_{t \in T} \sum_{h \in H} \sum_{g \in G} QH_{hgt} \times TH_{hg} + \sum_{t \in T} \sum_{g \in G} \sum_{o \in O} QG_{got} \times TM_{go}$$

$$+ \sum_{t \in T} \sum_{o \in O} \sum_{l \in L} QO_{olt} \times TO_{ol} + \sum_{t \in T} \sum_{g \in G} INV_{gt} \times HCG_g + \sum_{t \in T} \sum_{g \in G} \sum_{o \in O} QG_{got} \times (1 - \widetilde{\Theta_o}) * TC_o$$

*Subject to:*

$$\sum_{s \in S} QS_{swpt} \times (1 - \bar{p}) \geq \sum_{h \in H} QW_{whpt} + \sum_{g \in G} QD_{wgpt} \times \bar{p} \quad \forall p, w, t \quad (2)$$

$$\sum_{w \in W} QW_{whpt} \times (1 - \bar{p}) + QU_{pt} \geq \widetilde{D}_{hpt} \quad \forall p, h, t \quad (3)$$

$$\sum_{w \in W} QS_{swpt} \leq CS_{sp} \times S_{st} \quad \forall p, s, t \quad (4)$$

$$\sum_{h \in H} QW_{whpt} \leq CW_{wp} \times W_{wt} \quad \forall p, w, t \quad (5)$$

$$\sum_{w \in W} QD_{wgpt} \times \bar{p} + \sum_{h \in H} QH_{hgt} \leq CG_g \times G_{gt} \quad \forall p, g, t \quad (6)$$

$$\sum_{g \in G} QG_{got} \times (1 - \widetilde{\Theta}_o) \leq CO_o \times O_{ot} \quad \forall o, t \quad (7)$$

$$\sum_{o \in O} QO_{olt} \leq CL_l \times L_{lt} \quad \forall l, t \quad (8)$$

$$\sum_{g \in G} QH_{hgt} = \widetilde{VW}_{ht} \quad \forall h, t \quad (9)$$

$$INV_{gt} = INV_{g(t-1)} + \sum_{w \in W} \sum_{p \in P} QD_{wgpt} \times \bar{p} + \sum_{h \in H} QH_{hgt} - \sum_{o \in O} QG_{got} \times (1 - \widetilde{\Theta}) \quad \forall g, t \quad (10)$$

$$INV_{gt} \leq CG_g \times G_{gt} \quad \forall g, t \quad (11)$$

$$\sum_{h \in H} QH_{hgt} + \sum_{w \in W} \sum_{p \in P} QD_{wgpt} \times \bar{p} \leq \sum_{o \in O} QG_{got} \quad \forall g, t \quad (12)$$

$$\sum_{g \in G} QG_{got} \times (1 - \widetilde{\Theta}_o) \leq \sum_{l \in L} QO_{olt} \quad \forall o, t \quad (13)$$

$$\sum_{s \in S} QS_{swpt} \leq M \times Y_{wt} \quad \forall p, w, t \quad (14)$$

$$\sum_{h \in H} QH_{hgt} \leq M \times G_{gt} \quad \forall g, t \quad (15)$$

$$S_{st}, W_{wt}, G_{gt}, O_{ot}, L_{lt}, Y_{wt} \in \{0,1\} \; \forall \, t \in T, s \in S, w \in W, g \in G, o \in O, l \in L. \quad (16)$$

$$QS_{swpt}, QW_{whpt}, QD_{wgpt}, QU_{whpt}, QH_{hgt}, QG_{got}, QO_{olt}, INV_{gt} \geq 0 \quad (17)$$

Equation (1) minimizes the total cost of the SC; constraint (2) states that the order amount of vaccine items from the supplier should be greater than the amount of perfect and defective vaccine shipped from the Healthcare distribution centers to healthcare centers and from distribution centers to vaccine waste storages respectively. Constraint (3) maintain the demand in the healthcare center, constraint (4), (5),(6) ,(7), and (8) are capacity limitation for the supplier, warehouse vaccine waste storages, vaccine waste treatment centers and landfills respectively. (9) Ensure that the vaccine waste generated in the healthcare centers equal to the total amount of vaccine waste shipped to the storage center. (10) Represents the inventory balance at the vaccine waste storages,(11) limits the vaccine waste inventory to not exceed the storage center capacity, (12) is the flow balance of vaccine waste and defective vaccine items shipped from healthcare centers and distribution centers respectively to vaccine waste storages. (13) Is the flow balance of vaccine wastes between treatment centers and landfills. (14) States that if an order is placed, a cost is incurred for each supplier. (15) Ensure the assignment of healthcare centers to the opened storage centers. (16) Are binary constraints, (17) is non-negativity constraints.

It is assumed that all vaccine units received from the supplier(s) are screened and inspected in the distribution centers. However, the inspection process is not perfect, and two types of errors are committed: type I error ($\gamma_1$) is committed when misclassifying of a conforming vaccine unit occurs, and type II error ($\gamma_2$) is committed when a misclassifying of a non-conforming vaccine unit occurs. The probability of underlying errors is constant in all distribution centers. If the proportion of defective units in a vaccine order is denoted as p, then the apparent non-defective items probability is obtained as follows.

$$(1-p)(1-\tilde{\gamma}_1) + p\,\tilde{\gamma}_2 = 1 - \tilde{\gamma}_1 - p(1 - \tilde{\gamma}_1 - \tilde{\gamma}_2) = 1 - \bar{p} \qquad (18)$$

Where

$$\bar{p} = \tilde{\gamma}_1 + p(1 - \tilde{\gamma}_1 - \tilde{\gamma}_2) \qquad (19)$$

The apparent probability of defective items.

$$p(1-\tilde{\gamma}_2) + (1-p)(\tilde{\gamma}_1) = p - p\tilde{\gamma}_2 - p\tilde{\gamma}_1 + \tilde{\gamma}_1 = \tilde{\gamma}_1 + p(1 - \tilde{\gamma}_1 - \tilde{\gamma}_2) = \bar{p} \qquad (20)$$

## 3   ROBUST CONVEX OPTIMIZATION

Robust optimization (RO) is young research field which has been developed in the last 2 decades. RO is very practical, since it is tailored to the information at hand, and it leads to computationally tractable formulations, Gorissen et al.[8]. Different approaches have been discussed in literature; the most well-known concepts are adjustable robust optimization introduced by Ben-Tal et al.[9], while interval and budget uncertainty is coined by Bertsimas and Sim [10]. The proposed model is formulated following the approach provided by Li and al.[11], which has the flexibility of choosing the desired uncertainty (box, polyhedral, combined interval-polyhedral) set based on the decision maker's risk attitude.

### 3.1   Definition

By considering uncertainty in the deterministic model the robust counterparts can be formulated by introducing the uncertain parameter within the desired studied set (box, polyhedral, combined interval-polyhedral). In the following we show the combined interval-interval-polyhedral robust counterpart formulation based on Bertsimas and Sim [10], Li et al.[11]:

Let consider the following deterministic linear optimization problem:

$$max\ c^T x_j \qquad (21)$$

$$s.t.\ A\ x_j \leq b_j \qquad (22)$$

$$x \in X \qquad (23)$$

Where $c^T$, $A$, $b_j$ and $X$ are n-vector of coefficients, $m \times n$ matrix of coefficients consumptions, m-vector of the Right-Hand-Side (RHS) values, and the uncertainty set of solutions. Assume that $c^T$, $A(a_{ij})$ and $b_j$ are uncertain parameters and considered as bounded, symmetric and independent. Then these parameters take the values in $[c_j - \hat{c}_j, c_j + \hat{c}_j]$ , $[a_{ij} - \hat{a}_{ij}, a_{ij} + \hat{a}_{ij}]$, and $[b_j - \hat{b}_j, b_j + \hat{b}_j]$   where   $c_j$ ,$a_{ij}$ and $b_j$ represent the nominal value of the parameters; $\hat{c}_j$ and $\hat{a}_{ij}$ represent the constant perturbation.
Then the equation can be rewritten as follows:

$$\max \sum_{j \epsilon J} \tilde{c}_j \, x_j \quad (24)$$

$$s.t. \sum_{j \epsilon J} \tilde{a}_{ij} x_j \leq \tilde{b}_j \quad \forall \, i \epsilon \, I \quad (25)$$

Where $\tilde{c}_j$, $\tilde{a}_{ij}$ and $\tilde{b}_j$ represent the true values which are subject to uncertainty. Therefore, the above problem is formulated as follows:

$$\max \sum_{j \epsilon J} c_j \, x_j + z_0 \Gamma_0 + \sum_{j \epsilon J} p_{0j} \quad (26)$$

$$s.t. \sum_{j \epsilon J} a_{ij} x_j + [z_i \Gamma_i + \sum_{j \epsilon J} p_{ij} + p_{i0}] \leq b_j \quad \forall \, i \epsilon \, I \quad (27)$$

$$z_0 + p_{0j} \geq \hat{c}_j \, \forall \, j \epsilon \, J_0 \quad (28)$$
$$z_i + p_{ij} \geq \hat{a}_j \, \forall \, j \epsilon \, J_i \quad (29)$$
$$z_i + p_{io} \geq \hat{b}_j \, \forall \, j \epsilon \, J_i \quad (30)$$
$$z_0, z_i, p_{0j}, p_{ij}, p_{i0} \geq 0 \quad (31)$$

Where $z_0, z_i, p_{0j}, p_{ij}, p_{i0}$ are auxiliary variables, $J_i$ represents the index subset that contains the variable indices whose corresponding coefficients are subject to uncertainty. $\Gamma_i$ is the adjustable uncertainty set parameter for the combined interval-polyhedral uncertainty set.

### 3.2 Combined Interval-Polyhedral robust counterpart formulation

The proposed model is formulated as a robust convex optimization problem with three different uncertainty sets: box, polyhedral and combined interval-polyhedral respectively.

| Robust parameters | Robust variables |
|---|---|
| $\Gamma_{DEF}, \Gamma_{TR}, \Gamma_D, \Gamma_{VW}$ Polyhedral adjustable uncertain parameters for inspection errors, treatment rate, demand and vaccine waste respectively. <br> $\hat{p}, \widehat{\Theta}_o, \widehat{D_{hpt}}, \widehat{VW}_{ht}$ Parameter deviation from nominal values for inspection errors, treatment rate, demand and vaccine waste respectively. | $Z_7, Z_8, Z_9$ Combined interval-polyhedral auxiliary variables. <br> $m_{DEF}, W_{whpt}^{DEF}, W_{wgpt}^{DEF}, W_{swpt}^{DEF}$ Auxiliary variable to the change in inspection errors. <br> $m_{TR}, W_{got}^{TR}, W_{hgt}^{TR}$ Auxiliary variable to the change in treatment rate. <br> $m_D, W_{hpt}^D$ Auxiliary variable to the change in demand. <br> $m_{VW}, W_{ht}^{VW}$ Auxiliary variable to the change in vaccine waste. |

The combined interval-polyhedral robust counterpart is formulated by taking into consideration the uncertainties in the objective function and changing them into constraints; (for the proof see Li et al.[11]).

The constraints contain uncertain parameters are replaced by its robust counterparts then added to the deterministic model as follows.

$$Z_7 - \sum_{t \in T} \sum_{w \in W} \sum_{h \in H} \sum_{p \epsilon P} QW_{whpt} \times TW_{whp} \times (1 - \bar{p}) - \sum_{t \in T} \sum_{w \in W} \sum_{h \in H} \sum_{p \epsilon P} W_{whpt}^{DEF} - \Gamma_{DEF} \times m_{DEF} \geq 0 \quad (32)$$

$$Z_8 - \sum_{t \in T} \sum_{w \in W} \sum_{g \in G} \sum_{p \epsilon P} QD_{wgpt} \times TG_{wg} \times \bar{p} + \sum_{t \in T} \sum_{w \in W} \sum_{h \in H} \sum_{p \epsilon P} W_{wgpt}^{DEF} + \Gamma_{DEF} \times m_{DEF} \geq 0 \quad (33)$$

$$m_{DEF} + W_{whpt}^{DEF} \geq \hat{p} \times QW_{whpt} \quad \forall \, w, h, p, t \quad (34)$$
$$m_{DEF} + W_{wgpt}^{DEF} \geq \hat{p} \times QD_{wgpt} \quad \forall \, w, g, p, t \quad (35)$$
$$m_{DEF} \geq 0, W_{whpt}^{DEF}, W_{wgpt}^{DEF} \geq 0 \quad (36)$$

$$Z_9 - \sum_{t\in T}\sum_{g\in G}\sum_{o\in O} QG_{got} \times TC_o \times (1-\Theta_o) - \sum_{t\in T}\sum_{g\in G}\sum_{o\in O} W^{TR}_{got} - \Gamma_{TR} \times m_{TR} \geq 0 \quad (37)$$

$$m_{TR} + W^{TR}_{got} \geq \hat{\Theta}_o \times QG_{got} \quad \forall g,o,t \quad (38)$$

$$m_{TR} \geq 0, W^{TR}_{got} \geq 0 \quad (39)$$

$$\sum_{s\in S} QS_{swpt} \times (1-\bar{p}) - \sum_{s\in S} W^{DEF}_{swpt} - \Gamma_{DEF} \times m_{DEF} \geq \sum_{h\in H} QW_{whpt}$$
$$+ \sum_{g\in G} QD_{wgpt} \times \bar{p} + \sum_{g\in G} W^{DEF}_{wgpt} + \Gamma_{DEF} \times m_{DEF} \quad \forall w,p,t \quad (40)$$

$$m_{DEF} + W^{DEF}_{swpt} \geq \hat{p} \times QS_{swpt} \quad \forall s,w,p,t \quad (41)$$

$$m_{DEF} + W^{DEF}_{wgpt} \geq \hat{p} \times QD_{wgpt} \quad \forall w,g,p,t \quad (42)$$

$$m_{DEF} \geq 0, W^{DEF}_{swpt}, W^{DEF}_{wgpt} \geq 0 \quad (43)$$

$$\sum_{w\in W} QW_{whpt} \times (1-\bar{p}) - \Gamma_{DEF} \times m_{DEF} - \sum_{w\in W} W^{DEF}_{whpt} + W^D_{hpt} + \Gamma_D \times m_D + QU_{hpt}$$
$$\geq D_{hpt} \quad \forall h,p,t \quad (44)$$

$$m_{DEF} + W^{DEF}_{whpt} \geq \hat{p} \times QW_{whpt} \quad \forall w,h,p,t \quad (45)$$

$$m_D + W^D_{hpt} \geq \hat{D}_{hpt} \quad \forall h,p,t \quad (46)$$

$$m_{DEF}, m_D \geq 0, W^D_{hpt}, W^{DEF}_{whpt} \geq 0 \quad (47)$$

$$\sum_{w\in W} QD_{wgpt} \times \bar{p} + \Gamma_{DEF} \times m_{DEF} + \sum_{w\in W} W^{DEF}_{wgpt} + \sum_{h\in H} QH_{hgt} \leq CG_g \times G_{gt} \quad \forall g,p,t \quad (48)$$

$$m_{DEF} + W^{DEF}_{wgpt} \geq \hat{p} \times QD_{wgpt} \quad \forall w,g,p,t \quad (49)$$

$$m_{DEF} \geq 0, W^{DEF}_{wgpt} \geq 0 \quad (50)$$

$$\sum_{g\in G} QG_{got} \times (1-\Theta_o) - \Gamma_{TR} \times m_{TR} - \sum_{g\in G} W^{TR}_{got} \leq CO_o \times O_{ot} \quad \forall o,t \quad (51)$$

$$m_{TR} + W^{TR}_{got} \geq \hat{\Theta}_o \times QG_{got} \quad \forall g,o,t \quad (52)$$

$$m_{TR} \geq 0, W^{TR}_{got} \geq 0 \quad (53)$$

$$\sum_{g\in G} QH_{hgt} + W^{VW}_{ht} + \Gamma_{VW} \times m_{VW} = VW_{ht} \quad \forall h,t \quad (54)$$

$$m_{VW} + W^{VW}_{ht} \geq \widehat{VW_{ht}} \quad \forall h,t \quad (55)$$

$$m_{VW} \geq 0, W^{VW}_{ht} \geq 0 \quad (56)$$

$$INV_{gt} = INV_{g(t-1)} + \sum_{w\in W}\sum_{p\in P} QD_{wgpt} \times \bar{p} + \Gamma_{DEF} \times m_{DEF} + \sum_{w\in W} W^{DEF}_{wgpt}$$
$$+ \sum_{h\in H} QH_{hgt} - \sum_{o\in O} QG_{got} \times (1-\Theta_o) + \Gamma_{TR} \times m_{TR} + \sum_{o\in O} W^{TR}_{got} \quad \forall g,p,t \quad (57)$$

$$m_{DEF} + W^{DEF}_{wgpt} \geq \hat{p} \times QD_{wgpt} \quad \forall w,g,p,t \quad (58)$$

$$m_{TR} + W^{TR}_{got} \geq \hat{\Theta}_o \times QG_{got} \quad \forall g,o,t \quad (59)$$

$$m_{TR}, m_{DEF} \geq 0, W^{TR}_{got}, W^{DEF}_{wgpt} \geq 0 \quad (60)$$

$$\sum_{h\in H} QH_{hgt} + \sum_{w\in W}\sum_{p\in P} QD_{wgpt} \times \bar{p} + \Gamma_{DEF} \times m_{DEF} + \sum_{w\in W}\sum_{p\in P} W^{DEF}_{wgpt} \leq \sum_{o\in O} QG_{got} \quad \forall g,o,t \quad (61)$$

$$m_{DEF} + W^{DEF}_{wgpt} \geq \hat{p} \times QD_{wgpt} \quad \forall w,g,p,t \quad (62)$$

$$m_{DEF} \geq 0, W^{DEF}_{wgpt} \geq 0 \quad (63)$$

$$\sum_{g\in G} QG_{got} \times (1-\Theta_o) - \Gamma_{TR} \times m_{TR} - \sum_{g\in G} W^{TR}_{got} \leq \sum_{l\in L} QO_{olt} \quad \forall o,t \quad (64)$$

$$m_{TR} + W^{TR}_{got} \geq \hat{\Theta}_o \times QG_{got} \quad \forall g,o,t \quad (65)$$

$$m_{TR} \geq 0, W^{TR}_{got} \geq 0 \quad (66)$$

Constraints (4),(5),(8),(11),(14),(15),(16) and (17) remain the same as in the deterministic.

## 4 RESULTS AND ANALYSIS

The uncertain parameters we consider are type I and type II errors ($\gamma_1, \gamma_2$) with 10% deviation of coefficient, demand $D_{hpt}$ and generated amount of waste in the healthcare centers $VW_{ht}$ with 15 % perturbation from nominal values, and finally a 10 % of deviation of defective treatment rate of waste $\Theta$. Moreover, the uncertain parameters are assumed to be subject to bounded uncertainty. The three-uncertainty sets considered are, box, polyhedral and combined interval polyhedral, the three proposed robust approaches are carried out for the uncertain parameters separately, and all of them simultaneously afterwards. Table-1 depicts the size of the adopted example with small size inputs.

Table 1: The size of the illustrative example.

|  | No. of time periods | No. of Products (vaccines) | No. of Suppliers | No. of distribution centers | No. of Healthcare centers | No. of Vaccine waste storages | No. of Vaccine waste Treatment centers | No. of land fills |
|---|---|---|---|---|---|---|---|---|
| Size of the example | 12 | 3 | 3 | 4 | 3 | 4 | 3 | 3 |

In order to provide stakeholders with a meaningful interpretation on the three approaches applied to our model, a comparison is conducted with an extensive computation based on the decision maker's risk attitude to help make better decisions. The problem is coded in Python 3.11.4 [12] and solved with Gurobi 10.0 [13] on Laptop Core(TM) i7 with 16 GB, 3.00 GHz of RAM. Table-2 shows the perturbation of the four random parameters for each uncertainty set in the worst-case scenario i.e.( $\Psi$=1, $\Gamma_i = |J_i|$ ) with 1%, 5% and 10% from the nominal values, moreover the robust objective function along with the average unmet demand, average vaccine inventory waste and price of robustness are computed, optimality gaps and solution time are calculated.

Table 2: Deviation of the uncertain parameters in the worst-case scenario for the three uncertainty sets.

|  | Deviation of the uncertain parameters | Robust objective function | Average Unmet demand | Average Vaccine Inventory-waste | Price of robustness % | Optimality gap % | CPU time |
|---|---|---|---|---|---|---|---|
| Box | 1% | 6.02055 e+07 | 57.29904 | 665.2994 | 29.87 | 0.0034 | 1 |
|  | 5% | 6.25635 e+07 | 65.93386 | 724.4574 | 34.95 | 0.0099 | 1 |
|  | 10% | 6.55994 e+07 | 76.77277 | 669.2883 | 41.5 | 0.0003 | 2 |
| Polyhedral | 1% | 6.72895 e+07 | 79.94379 | 144.2217 | 45.15 | 0.0090 | 3 |
|  | 5% | 9.96763 e+07 | 217.3719 | 322.8301 | 115.01 | 0.01 | 4 |
|  | 10% | 1.47974 e+08 | 400.5709 | 631.1597 | 219.2 | 0.0093 | 3 |
| Polyhedral + Interval | 1% | 5.643 e+07* | 49.29007 | 126.8978 | 21.72* | 0.01 | 195 |
|  | 5% | 5.83877 e+07* | 57.3647 | 137.235 | 25.95* | 0.0058 | 125 |
|  | 10% | 6.08649 e+07* | 67.0765 | 148.1362 | 31.29* | 0.0093 | 135 |

It can be concluded that the combined interval-polyhedral uncertainty set is less conservative in the worst-case scenario and gives robust solutions with 21.72 %, 25.95 and 31.29 % worse than the deterministic model, for the three deviations 1%,5% and 10% respectively. While box and polyhedral uncertainty sets are providing more conservative robust solutions in the same studied scenario. In addition, the combined interval-polyhedral uncertainty sets are able to provide less shortages and vaccine inventory waste in the worst-case scenario compared with other sets. Figure-2 depicts the average amount of unmet demand and vaccine inventory waste in terms of the uncertainty sets with different level of robustness. it can be observed that the polyhedral uncertainty set provides much worse results in terms of the average amounts of shortages and vaccine inventory waste compared to the other investigated sets.

From this point of view and based on figure-3 it can be observed that the combined interval-polyhedral set is less conservative when the level of robustness is approximately $\Gamma_i > 0.8$. On the other hand, when the set size parameter is around $\Psi \leq 0.8$ the box uncertainty set is less conservative and could provide good solutions.

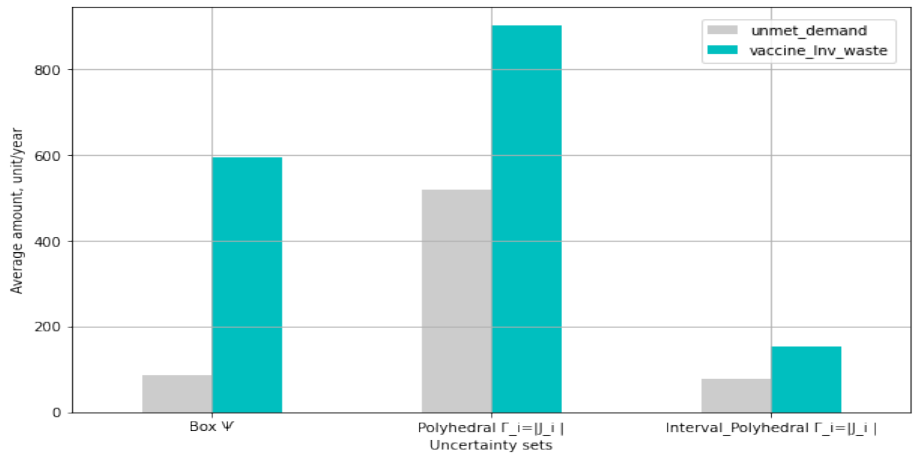

Figure 2: Average amount of unmet demand & vaccine inventory waste versus box, polyhedral, interval-polyhedral uncertainty sets.

In terms of solution time, the combined interval-polyhedral uncertainty set has much computational time compared to the other studied sets, while the optimality gap is showing marginal results for almost all the three sets.
Next, we study the combined interval-polyhedral set-in depth to see the effect of the four random parameters on the proposed model, and provide robust as well as feasible solutions for the problem.

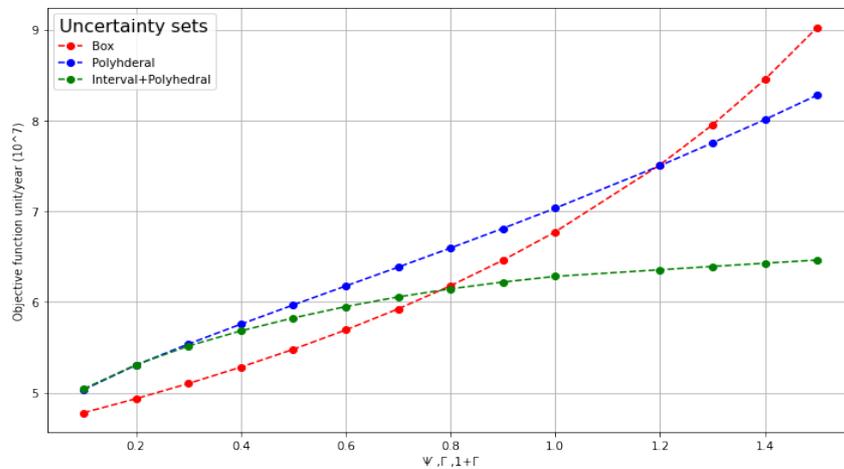

Figure 3: Total cost versus box, polyhedral, interval polyhedral with level of robustness.

Table 3: Combined Polyhedral-interval robust objective function for the four uncertain parameters (worst-case scenario $\Gamma_i = |J_i|$)

| Uncertain parameter | I and type II errors | Treatment rate | Demand | Vaccine Inventory waste |
|---|---|---|---|---|
| Objective function (Interval-polyhedral) | 5.67417 e+07* | 4.6556 e+07 | 4.9144 e+07 | 4.86644 e+07 |

The combined interval polyhedral uncertainty set approach is applied to the proposed model in the worst-case scenario for the uncertain parameters separately first, and then all simultaneously. From table-3, it can be concluded that type I&II uncertain parameters with $\Gamma_{DEF} = 1$ exert the biggest effect in terms of cost, with 22.4% worse than the deterministic, followed by demand, vaccine waste and defective treatment rate uncertainties with 6.01%, 4.97% and 0.42% of price of robustness respectively. To see the deterioration of the cost in terms the decision maker's risk attitude, type I&II and demand uncertainties are combined with defective treatment rate and vaccine waste random parameters. Figure-4&5 show the effect of incorporating the latter mentioned parameters with type I&II and demand uncertain parameters, we can observe the increase of the objective function gradually as $\Gamma_{DEF}$ and $\Gamma_D$ increase in the range $\Gamma_D \in [0, |J_i|]$ ,respectively.

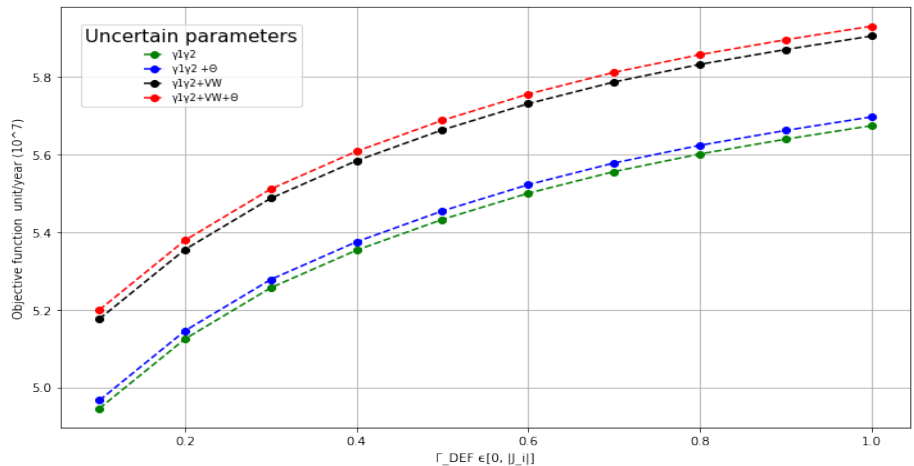

Figure 4: Objective function versus Γ_DEF (Interval Polyhedral uncertainty set).

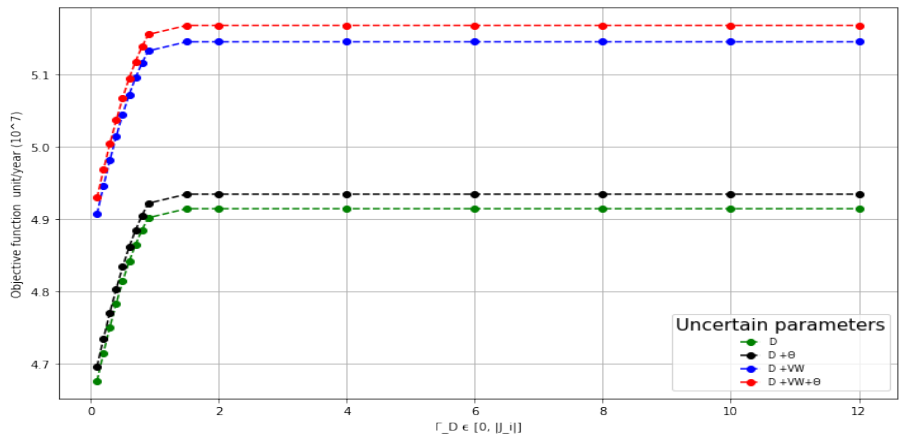

Figure 5: Objective function versus Γ_D (Interval-Polyhedral uncertainty set).

## 5 DISCUSSION

In this research, a multi-echelon multi-vaccine reverse supply chain problem is studied. In addition, a robust optimization approach is adopted to provide robust and feasible solution and protect against parameters deviation. First, we address all three uncertainty sets for our proposed problem, computational results show that the Interval-polyhedral uncertainty set is less conservative when $\Gamma_i > 0.8$ compared to other sets. However, when $\Psi_i \leq 0.8$ the box uncertainty is less conservative than polyhedral and combined interval-polyhedral uncertainty sets. Moreover, the three uncertainty sets are compared in terms of the average unmet vaccine demand and vaccine inventory waste, as expected the combined interval-polyhedral uncertainty set provides the best results with less amount of vaccine shortages and vaccine inventory waste among the three sets in the worst case-scenario. While the box uncertainty set with $\Psi_i = 0.5$ shows better solutions than others. The four uncertain parameters are further investigated based on combined interval-polyhedral uncertainty set. Extensive computations are conducted for the proposed approach. The computational results are able to provide robust solution for all uncertain parameters separately and simultaneously afterwards. Moreover, considering only deviation of type I&II parameters, with $\Gamma_{\text{DEF}} = |T|$ the solution is 22.4% worse than the deterministic, while demand uncertainty solution is 6.01% worse.

The obtained result also suggests that vaccine waste randomness has more pronounced effect on the objective function with 4.97% worse than the deterministic compared to defective treatment rate deviation, which has 0.42%.

The above results show the sensitivity of the proposed model (objective function) to parameters deviation, specifically to inspection errors randomness incorporated in the problem. In addition, taking into account uncertain parameters gradually deteriorate the total cost, while ignoring it could lead to sub-optimality or infeasibility.

It could be concluded that the proposed model is greatly affected by the uncertainty in inspection (type I&II errors), which calls for decision maker to be robust against uncertainty, shortages and disruptions in vaccine supply. Moreover, considering multiple uncertainties in the proposed model could improve the robustness of the supply chain network and help policymakers for better planning.

## 6 CONCLUSION

This manuscript contributes to the existing literature by proposing a new MILP model for multi-vaccine multi-echelon supply chain to help provide access to different types of vaccines during an outbreak such as COVID-19 and get dispose of infectious waste generated from massive vaccination campaign efficiently. While taking into consideration the decision maker's risk attitude by providing robust and feasible solutions in acceptable computational time in face of multiple uncertainties, such as inspection errors, demand, vaccine waste and defective treatment rate, which lie within box, polyhedral and combined interval-polyhedral uncertainty sets.

Further research could be focusing on optimizing the set size of the desired studied uncertainty set to provide less conservative solutions and account for the decision maker's risk behavior. Another interesting aspect of extending the on-hand problem by transforming the single objective into bi-objective or multi-objective function. While computation complexity could be handled by developing an exact or metaheuristic algorithm for a large-scale problem. Perishability of vaccines and its deterioration could be incorporated in the model along with inventory in healthcare centers.

**Declaration of Competing Interest**

The authors declare no conflict of interest.

**Acknowledgments**

The authors would like to acknowledge the support provided by the Deanship of Research Oversight and Coordination at King Fahd University of Petroleum & Minerals (KFUPM).